Optimal Regulator2 10 24 02

# Design of Optimal Regulators*

**Alexander A. Bolonkin**
NRC Senior Research Associate
Air Force Research Laboratory
1310 Avenue R, #6-F,
Brooklyn, NY 11229, USA
T/F 718-339-4563, aBolonkin@juno.com,
abolonkin@gmail.com, http://Bolonkin.narod.ru

**Robert L. Sierakowski**
Chief Scientist
Air Force Research Laboratory
AFRL/MNG 101 W. Eglin BLVD, Ste.105 ,
Eglin AFB, Florida 32542-6810

## Abstract

Current research suggests the use of a liner quadratic performance index for optimal control of regulators in various applications. Some examples include correcting the trajectory of rocket and air vehicles, vibration suppression of flexible structures, and airplane stability. In all these cases, the focus is in suppressing/decreasing system deviations rapidly. However, if one compares the Linear Quadratic Regulator (LQR) solution with optimal solutions (minimum time), it is seen that the LQR solution is less than optimal in some cases indeed (3-6) times that obtained using a minimum time solution. Moreover, the LQR solution is sometimes unacceptable in practice due to the fact that values of control extend beyond admissible limits and thus the designer must choose coefficients in the linear quadratic form, which are unknown.

The authors suggest methods which allow finding a quasi-optimal LQR solution with bounded control which is closed to the minimum time solution. They also remand the process of the minimum time decision.

----------------------------------

## Introduction

The LQR solution is easily and conveniently written using the Riccati equation as an optimal solution. The scientist who accepts this may be acting as an intoxicated man in a Russian anecdote: one night a man is observed creeping around a streetlight. A passerby asks him, what are you doing? – I lost money. Where did you lose the money? –There at the other end of the street. Then, why are you looking here? – This is where the light is!

The minimum time solution is more complex, however, it can be conveniently determined in many problems by the availability, generally, of high-speed computers. Also, this approaches us with a true minimum time solution.

For an $n$-dimensional problem with one control this solution found in general form in reference [1]. For the two-dimensional case this solution can be presented graphically, see ref. [1]. Methods for other general optimal solutions are offered in [2]-[4].

The LQR solution has three main issues:

1) The selection of the matrix coefficients in the performance index are designer selected and the solution is dependent upon the value of these coefficients.
2) The range of control values can be large in number and this not admissible for practice.

3) The “optimal” LQR solution can be up to 3-6 times worse, then the minimum time solution (see the example in this paper).

If a researcher chosees to use the LQR solution, the authors suggest a method for limiting maximum control (see point 2) as well as for the choice of selecting the coefficients in the performance index. This allows up to a 2-3 times improvement in the performance index (see accompanying examples) and thus makes the LQR solution acceptable in practical applications.

The traditional approach used in the design of a controlled structural system is to design the structure first by satisfying given requirements and then to design the control system. The structure is designed with such constraints placed on weight, allowable stresses, displacements, buckling, general instability, frequency distributions, etc. When the selection of the geometry, cross-sectional area of the members, and material are determined for a specified structure, then the structural frequencies and vibration modes become important input in the design of the control system. Some investigators have written papers discussing an integrated design approach for optimal control. In most references, the control design procedures used, do not take into consideration the limitations on the control forces developed by the actuators, and have not been treated as constraints or design variables. In this paper the problems associated with the selection of the performance index, parameters, weight coefficient in the LQR problem, and limitation of control forces are addressed.

In the following sections, theories for the synthesis of an optimal control laws with a quadratic performance index and bounded control forces are given. This is followed by a SISO (Single Input, Single Output) control problem designed using both approaches for comparison of the end state trajectories, with different bounds placed on control forces. Next, the control system for an idealized wing-box is used to illustrate a design application of the method. A discussion on the application of a control system with bounded control for an integrated design of a structure and control system can be found in ref. [5]. Related articles are [6]-[10].

The article with all figures has size 5.8 Mb. The Arxiv gives only the space 1 Mb. That way the part of figures are deleted.

## 1. Optimal Control

The general optimal control problem can be described by the following equations

$$I = F_0(x_1, x_2) + \int_{t_1}^{t_2} f_0(t,x,v)dt \ , \quad dx/dt = f(t,x,v), \quad x(t_1) = x_1, \quad x(t_2) = x_2 \qquad (1\text{-}1)$$

where $I$ is the functional (objective function), $t$ is time, $x$ is a $n$-dimensional vector of state, and $v$ is a $p$-dimensional vector of control forces. The vector $v \in V$ where $V$ can be a bounded domain. Boundary conditions $t_1, t_2, x_1, x_2$ are usually given, $(t_1, t_2) \in T$.

The control parameter, $v$ is calculated so that $I = min$. To find the solution to this problem, assume the function

$$\psi = \psi(t,x) \qquad (1\text{-}2)$$

and write the new functional

$$J = A_1 + \int B_1 dt \qquad t \subset [t_1, t_2], \qquad (1\text{-}3)$$

where

$$A_1 = F_o + \psi(t_2) - \psi(t_1)\ , \qquad B_1 = f_o - (\partial\psi/\partial x)f - (\partial\psi/\partial t)\ . \qquad (1\text{-}4)$$

Here $(\partial\psi/\partial x)$ is $n$-dimensional vector of partial derivatives. The global minimum is

$$A_0 = \inf_{x_1,x} A_1(x_1,x_2), \quad B_0 = \inf_{v,x} B_1(t,x,v) \quad for \forall t \in T\ . \qquad (1\text{-}5)$$

Depending on the nature of the functions used for $\psi$, a different set of algorithms for obtaining the infimum can be developed. For example, if Eq.(1-2) takes the form

$$\psi = \lambda(t)x\ , \qquad (1\text{-}6)$$

where $\lambda(t)$ is an $n$-dimensional vector, the global minimum functions can be written as,

$$A_o = inf\ (F_1 + \psi(t_2) - \psi(t_1)), \qquad B_o = inf\,[f_o - \lambda f(t,x,v) - (d\lambda/dt)x] = inf\ B\ . \qquad (1\text{-}7)$$

$x_1,x_2$ $x,v$ $x,v$

Using $\partial B/\partial x =$ 0 and Eq.(1-7) gives

$$d\lambda/dt = -\,\partial H/\partial x, \quad B = \inf B, \quad v \in V, \qquad (1\text{-}8)$$

where $H = \lambda f(t,x,u) - f_o$ .

Eq. (1-8) can be integrated to find $\psi$, to obtain the optimal control $v$ and the optimal trajectory $x(t)$. Another way is to enforce the condition

$$B = \inf_{v \in V}\left[ f_0 - \frac{\partial \psi}{\partial x} f - \frac{\partial \psi}{\partial t} \right] = 0 , \qquad (1\text{-}9)$$

everywhere in the admissible domain for $x$. In this case, the equation for particular derivatives can be solved and the syntheses of the optimal control $v = v(t,x)$ and the field of the optimal trajectory in the admissible domain is obtained.

The two control design approaches with constraints on the maximum control forces are discussed in this section. In the first section an objective function for establishing of the minimum time to suppress vibration is discussed and in the second, the quadratic function is minimized.

### A. Minimum Time

Since the main purpose of the controller is to suppress vibrations in minimum time, the time for the system to come to rest is taken as the objective function. A functional expression for this can be written

$$I = \int_0^T dt, \qquad T = min \qquad (1\text{-}10)$$

subject to

$$dx/dt = Ax + bf, \quad x(0) = x_o, \quad x(T) = 0 \qquad (1\text{-}11)$$

with control force limits

$$|f_i| \leq F_i \;, \quad i = 1,2,...,p \;. \qquad (1\text{-}12)$$

This problem can be written in short form as

$$\min I = \int_0^T dt, \qquad dx/dt = Ax + bf, \; x(0) = x_o, \; x(T) = 0, \; |f| \leq F, \qquad (1\text{-}13)$$

where $x$ is the state vector of dimension $2n$. $A$ is the $2n\times 2n$ plant matrix, $B$ is $2n\times p$ control matrix, $f$ is the control force vector of dimension $p$, $x(0)$ is the initial state vector, and $x(T) = 0$ is the final state of the system. $B_o$, in Eq.(1-7), for this problem can be written as

$$B_o = 1 - \Sigma(\partial\psi/\partial x_i)(dx_i/dt) - (\partial\psi/\partial t) \quad (i = 1.,,, n). \qquad (1\text{-}14)$$

Substituting

$$\psi = \Sigma\lambda_i(t)x_i \qquad (i = 1,2,...,2n) \qquad (1\text{-}15)$$

and Eq.(1-11) into Eq. (1-14) gives

$$B = 1 - \sum_{J=1}^{2n} \lambda_J \left( \sum_{i=1}^{2n} a_{ij} + \sum_{k=1}^{p} b_{jk} f_k \right) - \sum_{j=1}^{2n} \dot{\lambda}_j x_j \quad . \qquad (1\text{-}16)$$

Taking the partial derivatives of $B$ ($\partial B/\partial x_i$) gives

$$d\lambda_j/dt = -\sum_j a_{ij}\lambda_j \qquad i = 1,2,\ldots,n \,, \; j = 1, 2,\ldots,2n \;. \qquad (1\text{-}17)$$

$min\ B_o$ gives the control force $f$

$$f_i = |F_i| \,\mathrm{sign}\, (\Sigma b_{jk}\lambda_k) \qquad i = 1, 2,\ldots,p \,; \; k = 1, 2,\ldots,2n \,. \qquad (1\text{-}18)$$

Using Eqs.(1-13),(1-17) and (1-18), the optimal control force $f_i(t)$ and trajectory $x_i(t)$ can be calculated. However the initial $\lambda_i(0)$ for our trajectory with $x(0) = x_o$ is not known. To find $\lambda_i(0)$, any suitable gradient method can be used. For example, if the assume some initial state $\lambda_i(0)$ and integrate Eqs.(1-13),(1-17) and (18), we can calculate the function

$$I = T + \sum_{i=1}^{2n} C_i x_i^2(T), \qquad C_i > 0 \ , \qquad (i=1, 2, \ldots, 2n). \qquad (1\text{-}19)$$

Here $C_i$ are weight coefficients. If $\sum C_i x^2{}_i(T) < C_0$, where $C_o$ is small, the problem can be considered as solved. Time is optimal and $x_i(t)$ is the optimal trajectory which satisfies the final condition $x_i(T) = 0$. If $\sum C_i x^2{}_i(T) > C_o$ we can choose a new $\lambda i(0)$ by any method and repeat the process until it satisfies $\sum C_i x^2{}_i(T) < C_o$.

In practice, a new independent variable $\tau$ is introduced as $t = c\tau$, which can be included with Eq.(1-11) to prouder the additional equation

$$dt/d\tau = c \ . \qquad (1\text{-}20)$$

Additionally, introducing a fixed interval of integration $[0, \tau_1]$ a new set of equations become

$$\min I = \int_0^{\tau} c d\tau \ , \quad dx/dt = (Ax + Bf)c \ . \ \ x(o) = x_o, \ \ x(\tau_1) = 0, \ \ |f| \leq F \ , \qquad (1\text{-}21)$$

where $c$ is some constant, which is also selected. Eq.(1-19 ) thus becomes

$$I_1 = \Sigma C_i x^2{}_i(\tau_1) \ , \quad i = 1,2,\ldots,2n \quad . \qquad (1\text{-}22)$$

For the structural system as defined by Eqs.(1-11)-(1-12) this problem can be solved for the case in which the number of control inputs, $p$, is equal to the number of modeled structural degrees of freedom, $n$. However, numerical difficulties would be encountered when this condition is not satisfied. Typical difficulties would be the occurrence of many local minimums, poor convergence, and the need for smaller step sizes.

### B. Linear quadratic regulator (LQR) with bounded control

In this case, a performance index, $J$, is defined as

$$J = \int_0^{\infty} (x^T \underline{Q} x + f^T \underline{R} f) dt \qquad t \in [0, \infty] \ . \qquad (1\text{-}23)$$

Where $\underline{Q}$ and $\underline{R}$ are state and control weighting matrices. The matrix $\underline{Q}$ must be positive semi definite ($x^T\underline{Q}x \geq 0$), and $\underline{R}$ must be positive definite ($f^T\underline{R}f > 0$). The dimensions of $\underline{Q}$ and $\underline{R}$ depend on the size of vectors the $x$ and $f$, respectively. The matrices $\underline{Q}$ and $\underline{R}$ can be written as

$$\underline{Q} = \delta Q \qquad (1\text{-}24)$$

and

$$\underline{R} = (1/\gamma) R^{-1} \qquad (1\text{-}25)$$

where $\sigma$ and $\gamma$ are the design positive variables and $\underline{Q}$ and $R^{-1}$ are constant identity matrices. The weighting matrix $R$ is defined in terms of the inverse of the constant matrix $R$ in order to maintain positive definiteness. The function $B$, Eq. (1-9) for the performance index defined in Eq.(1-23) and the constraint equation Eq.(1-11), become

$$\overline{B} = \inf_f \left[ (\delta x^T Qx + \gamma f^T Rf) - \frac{\partial \psi}{\partial x}(Ax + Bf) - \frac{\partial \psi}{\partial t} \right] = 0 \qquad (1\text{-}26)$$

If V is an open domain, the function $\psi$ can bewritten

If $V$ represents an open domain, the function $\psi$, can be written in the form

$$\psi = x^T P x \ , \qquad (1\text{-}27)$$

where $P$ is a $2n$-dimensional unknown matrix.

Substituting $\psi$ Eq. (1-27) into Eq, (1-26), we obtain the equation

$$\sigma\gamma \underline{Q} + PA + A^T P - \gamma PBRB^T P = 0 \quad . \qquad (1\text{-}28)$$

Equation (1-28) is the Riccati equation. A solution of this equation gives the matrix $P$ and one can find the optimal control force as

$$f = - Gx \qquad (1\text{-}29)$$

where

$$G = \gamma R B^T P \qquad (1\text{-}30)$$

Integrating Eq. (1-11) using Eq. (29) to obtain the optimal trajectory for the LQR functional. Eq. (1-29) may give unrealistic values of control depending on the selection of γ. The magnitude of control can be decreased by increasing γ, however, this may cause other perturbations of the system (such as the time it takes the oscillation to decay) to deteriorate.

In order to obtain more realistic results, bounds can be placed on the control force. This can be written as

$$|f_i| \le F_i, \quad F_i = \text{const}, \; i = 1, 2, \ldots, p \qquad (1\text{-}31)$$

where $F_i$ is the magnitude bounding each controller. To obtain an optimal solution, the following restrictions must be satisfied: (1) among these optimal synthesis of the control must exist in the domain of interest, (2) the function $B$ Eq. (1-30) must be convex, and (3) the limits of $F$ may be constant or dependent on time only and $F$ must not be equal to zero at any time (Note: if $F$ is very small a loss in stability can occur). For a solution, the system of Eqs. (1-11) and (1-29) must be integrated along with limits imposed by equation (1-31).

The norm for the displacements or total deviation can be defined by

$$R_x(t) = S = [\sum_{i=1}^{n} x_i^2(t)]^{1/2} \qquad i = 1, 2, \ldots, n \qquad (1\text{-}32)$$

This norm is zero at the time the deviation is zero, and the structure stops vibrating. In the LQR solution domain this time equals infinity. For studying the behavior and comparison of different control systems, a measure of performances has been used based upon. The time required to reduce the norm of the displacements to 2% of their initial value.

## Numerical Examples.

**Example 1. SISO problem.**

For comparison of systems with different objective functions, a vibrating structure with a single physical degree of freedom was been investigated. This system is described by equation the following set of

$$dx_1/dt = x_2, \quad dx_2/dt = -\omega^2 x_1 - 2\zeta\omega x_2 + cf, \; x_1(0) = 0, \; x_2(0) = 1, \; |f| \le 1 \qquad (1\text{-}33)$$

where $\omega = 2$ is the frequency, $\zeta = 0.03$ is the damping, $c = 1$, and $|f| \le 1$ is the control.

The problem is solved having an objective function for minimum time as

$$\min T = \int_0^T dt, \; x_1(T) = 0, \; x_2(T) = 0 \qquad (1\text{-}34)$$

Eqs. (1-17) and (1-18) for the system defined in Eq. (1-33) become

$$d\lambda_1/dt = -\omega^2\lambda_2, \quad d\lambda_2/dt = \lambda_1 - 2\xi\omega\lambda_2, \quad f = |F| \operatorname{sign} \lambda_2 . \qquad (1\text{-}35)$$

Eqs. (1-33)-(1-35) are integrated and the initial values $\lambda_1(0)$, $\lambda_2(0)$ are chosen such that the conditions $x_1(T) = x_2(T) = 0$ are satisfied. The details of the solution scheme are not given here because of space limitations.

The performance for the linear quadratic regulator (LQR) is

$$J = \int_0^\infty \frac{1}{2}\left[(\delta_1 x_1^2 + \delta_2 x_2^2\right]dt \; . \qquad (1\text{-}36)$$

Using this performance index and solving the Riccati Eq. (1-29) gives

$$f = 2(c/\gamma)(c_{12}x_1 + c_2 x_2), \qquad (1\text{-}37)$$

where

$$c_{12} = -[\omega^2 + (\omega^4 + c_o\delta_1)^{0.5}]/2c_o, \; c^2 = \{-\zeta\omega + [\xi^2\omega^2 + (0.25\delta_2 + c_{12})c_o]^{0.5}\}/c_o, \quad c_o = c^2/\gamma .$$

In the case of $\delta = \delta_1 = \delta_2$, the time history depends only on $\delta/\gamma$. The total deviation is

$$R_x = S = (x_1^2 + x_2^2)^{1/2} . \qquad (1\text{-}38)$$

Eq. (1-33) is integrated with control given in Eq. (1-37).

The results of this investigation for the case $T$ = min, $\sigma/\gamma$ = 0.25 and 100 and no control (open-loop system) are shown in Fig's 1, 2, & 3.

Fig. 1 shows the time history of deviation of $x_2$. As can be seen, an LQR with $\sigma/\gamma$ = 100 gives better results ($t$ = 4 sec) than an LQR with $\sigma/\gamma$ = 0.25 (time is more than 15 sec) however an even better result is obtained with an objective function of minimum time. In the last case, oscillations are terminated in 1.5 sec.

Figure 2 shows the variation of a bounded control force $|f| \leq 1$ for the case of $T$=min, LQR when $\sigma/\gamma$ = 0.25 and $\sigma/\gamma$ = 100. The case LQR ($\sigma/\gamma$ = 0.25) does not use the full control force, the case LQR ($\sigma/\gamma$ =100) uses more of the control force, and case $t$ =min uses the maximum control force all the time.

Fig.3 shows the time history for the total deviation ($R_x$) with no control, with an objective function for minimum time and with LQR given by control bounds $|f| \leq 1$.

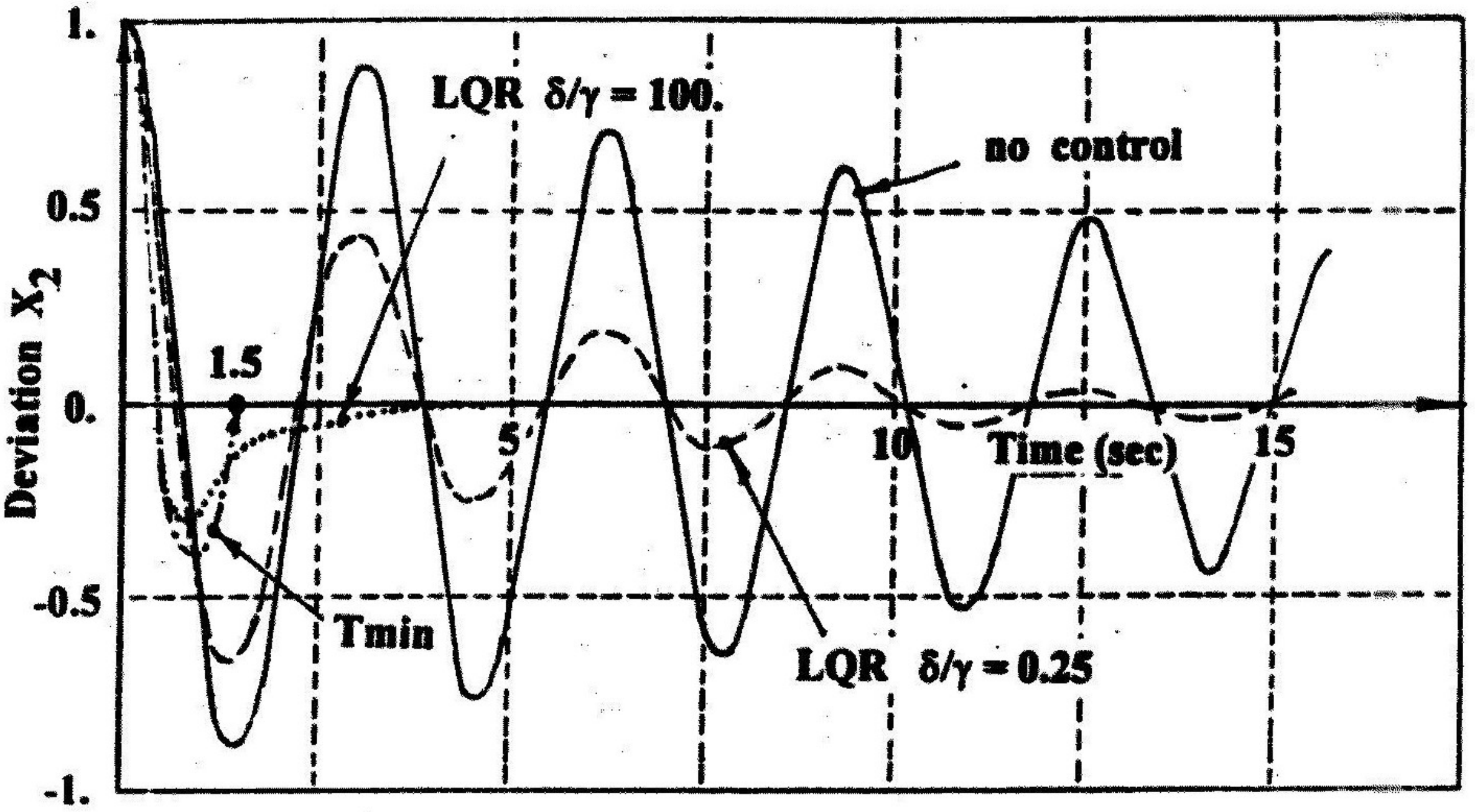

Fig 1

A structural system with any number of degrees of freedom can be transformed into pairs of equations (1-33)(see later Eq.(1-40)-(1-48)) where every pair is independent from the other. If the number of controls equals the number of degrees of freedom the design approach based on minimum time can be used. However, if the number of controls is less than the number of pairs of equations, the solution for the functional $T$ = min becomes very complex. In this case, the LQR approach is a variable alternative.

Fig.2 (deleted)

Fig.3 (deleted)

**Example 2. Wing Box**

In order to illustrate the application of an approach using the linear quadratic regulator with bounded control, the wing box problem in reference [5] is used and shown in Fig. 4. This structure has thirty-two elements and twenty-four degrees of freedom. The structure is a cantilever wing box idealized with bar elements capable of carrying axial loads only.

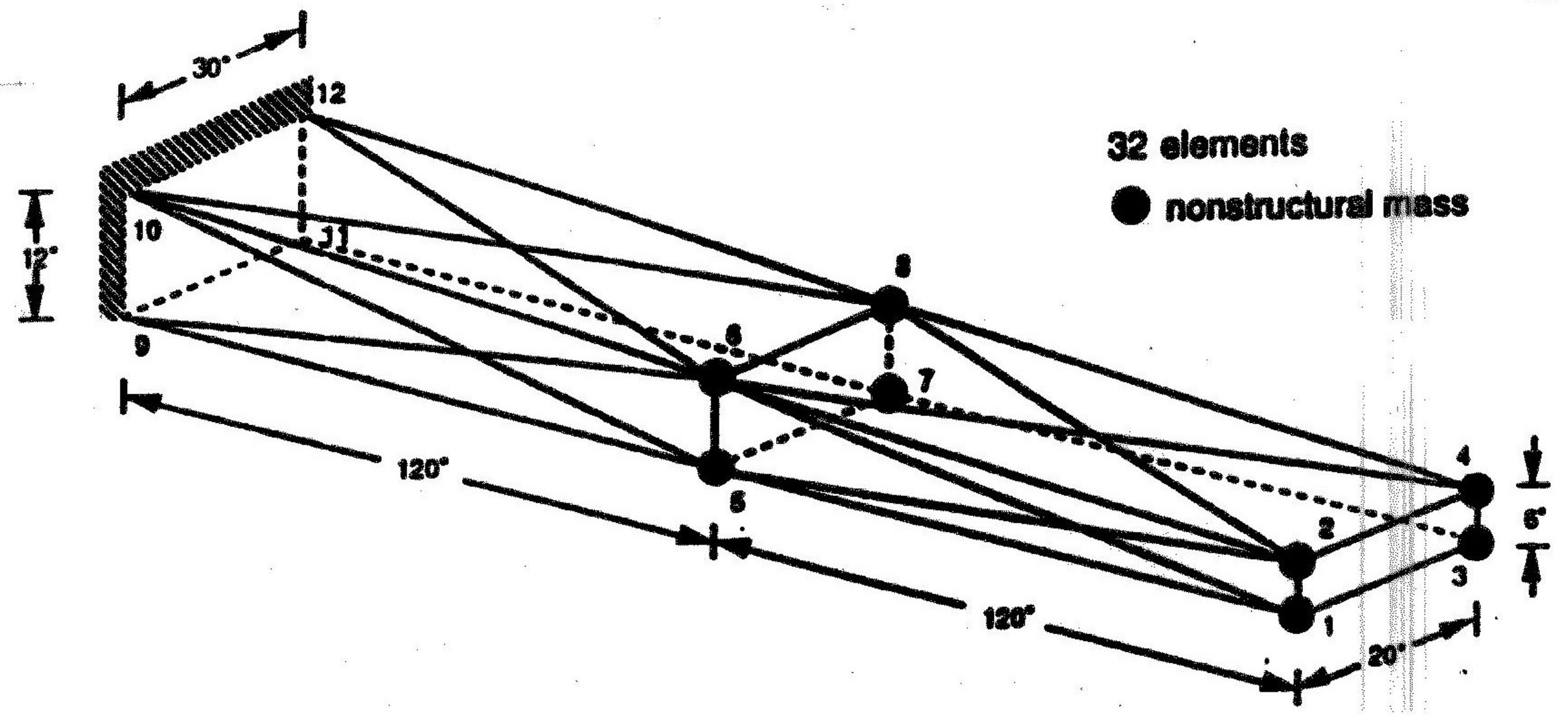

Fig.4.

The equations of motion for a flexible structure with no external disturbance can be written as

$$M\ddot{u} + E\dot{u} - Ku = Df\,, \qquad (1\text{-}39)$$

where $M$ is the mass matrix, $E$ is the damping matrix, and $K$ is the total stiffness matrix. These matrices are $n_1 \times n_1$, where $n_1$ is the number of degrees of freedom of the structure. In Eq. (1-39), $D$ is the applied load distribution matrix relating the control input vector $f$ to the coordinate system. The number of elements in vector $f$ is equal to the number of actuators, $p$. The vector u in Eq.(1-39) defines the structural response.

The coordinate transformation

$$u = [\Phi]\eta \qquad (1\text{-}40)$$

is introduced where η is the modal coordinate system and [Φ] is the $n_1 \times n_1$ modal matrix. Using Eq. (1-40), Eq.(1-39) can be transformed into $n_1$ uncoupled equations. These can be written as

$$\overline{M}\ddot{\eta} + \overline{E}\dot{\eta} + \overline{K}\eta = [\Phi]^T Df \qquad (1\text{-}41)$$

where

$$\begin{aligned} \overline{M} &= I = [\Phi]^T M[\Phi] \\ \overline{E} &= [2\zeta\omega] = [\Phi]^T E[\Phi] \qquad (1\text{-}42,\ 1\text{-}44) \\ \overline{K} &= [\omega^2] = [\Phi]^T K[\Phi] \end{aligned}$$

The matrices $\overline{M}, \overline{E},$ and $\overline{K}$ are diagonal square matrices, ω is the vector of structural frequencies, and ζ is the vector of modal damping factors. The modal matrix [Φ] is normalized with respect to the mass matrix. The control analysis is performed by reducing the second-order uncoupled equation [Eq.(1-41)] to a first-order equation. Only $n$ of $n_1$ uncoupled equations are used for the control system design. This can be achieved by using the transformation

$$x_{2n} = \begin{bmatrix} \eta \\ - \\ \dot{\eta} \end{bmatrix}_{2n} \qquad (1\text{-}45)$$

where $x$ is the state variable vector of size $2n$. This gives

$$\dot{x} = Ax + Bf \qquad (1\text{-}46)$$

where $A$ is a $2n \times 2n$ matrix and $B$ is a $2n \times p$ input matrix. The A matrix and the input matrices are given by

$$A = \begin{bmatrix} 0 & I \\ -\omega^2 & -2\xi\omega \end{bmatrix} \tag{1-47}$$

$$B = \begin{bmatrix} 0 \\ \Phi^T D \end{bmatrix} \tag{1-48}$$

The state output equation is given by

$$y = Cx \tag{1-49}$$

where $y$ is a $q \times 1$ output vector, $C$ is a$q \times 2n$ output matrix, and $q$ is equal to the total number of sensors. If the number of sensors and actuators equal and collocated, then $q = p$ and

$$C = B^T . \tag{1-50}$$

For this structure, Young's modulus and weight density are assumed to be equal to 10.5 x $10^6$ lbs/in$^2$ and 0.1 lbs/in$^3$, respectively. The actuators and sensors are assumed embedded in the structural elements and are collocated. The actuators are assumed to apply forces along the axial directions providing both out of plane, in plane and twist control for the structure. It is assumed that all structural modes have 1% structural damping and thus $\zeta$ in Eq. (1-9) was 0.01.

The control system utilizes four actuators and sensors collocated in the four members at the tip of the structure connecting nodes 1-2, 3-4, 1-3 and 2-4 respectively. Non-structural masses are located at nodes 1 through 8. Their magnitudes are 0.5 slugs at nodes 1 and 2; 1.5 slugs at nodes 3 and 4; 2.5 slugs at node 5 and 7.0 and 1.0 slugs at nodes 6 and 8 respectively. For the 24 structural degrees of freedom, the full order state space matrix in Eq. (1-11) is 48 x 48. Since there are four actuators and sensors, the input matrix $B$ and output matrix $C$ are 48 x 4 and 4x48, respectively. The cross-sectional areas of the rod elements were equal to 0.1 in$^2$. The weighting matrices $Q$ and $R$ in Eq. (1-28), (1-29) were equal to the identity matrix.

The four values of the weighting parameter ratios $\sigma/\gamma$ selected for this study are 0.1, 1.0, 100 and 1000, respectively. The maximum control forces generated by the four actuators are given in Table 1.

**Table 1**. Calculated cases

| Control bound $F$ | $\infty$ | 0.5 | 0.15 | 0.05 |
|---|---|---|---|---|
| $\delta/\gamma = 0.1$ | + | + | + | + |
| $\delta/\gamma = 1$ | + | + | + | + |
| $\delta/\gamma = 100$ | + | + | + | + |
| $\delta/\gamma = 1000$ | + | + | + | + |
| No control | + | | | |

The initial condition used for designing the controllers is a unit displacement at node 1 in the z-direction. This condition is used for all cases and also to obtain the response curves. The response curves are given for only a few cases because of space limitations. The three limits on the maximum allowable control forces are set equal to 0.5, 0.15, and 0.05 respectively. The different cases considered are summarized in Table 2.

**Table 2**. Maximum actuator forces

| | Actuator # | | | |
|---|---|---|---|---|
| Value $\delta/\gamma$ | 1 | 2 | 3 | 4 |
| $\delta/\gamma = 0.1$ | 0.05 | 0.05 | 0.07 | 0.03 |
| $\delta/\gamma = 1.0$ | 0.20 | 0.24 | 0.31 | 0.12 |
| $\delta/\gamma = 100$ | 1.31 | 1.89 | 3.30 | 1.23 |
| $\delta/\gamma = 1000$ | 2.95 | 2.25 | 8.25 | 3.73 |

In the case of $\sigma/\gamma = 0.1$ the maximum actuator forces are less than 0.15, and for $\sigma/\gamma = 1.0$, they are less than 0.5. Fig. 5 shows the time history of the displacement norm without control bound for the four values of $\sigma/\gamma$ and without control.

Fig.5 (deleted)

The maximum value of the displacement norm as a function of time is shown in Fig. 6.
The time required to decrease the displacement norm to 2% of its initial value 1.0 is shown in Fig. 7.

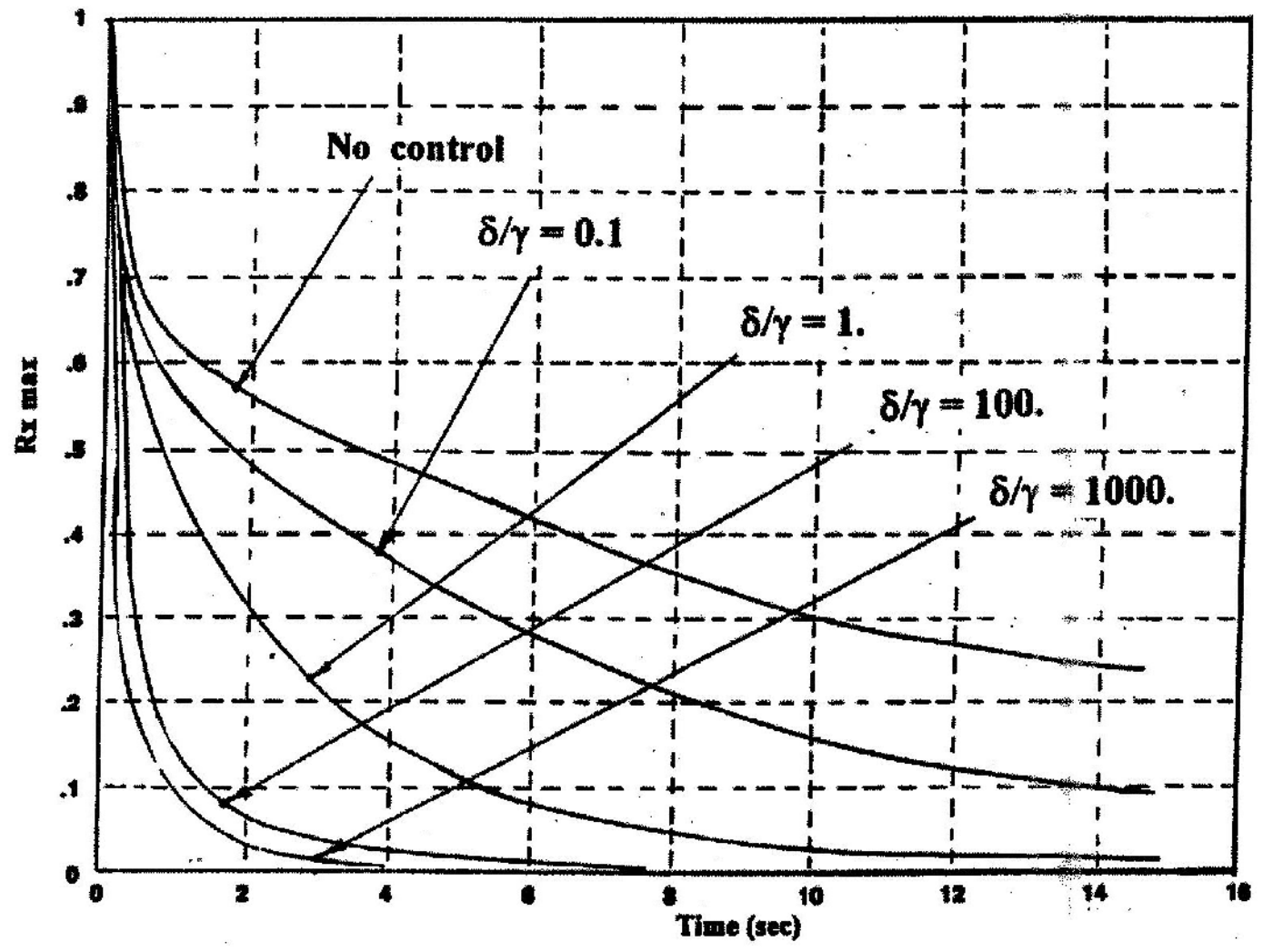

Fig.6. Fig.7 (deleted)

In the case of no-control, the total time needed to reduce the displacement norm to two percent of the initial value is larger than 100 seconds. The variation in the control force in actuator 1 as a function of time for $\sigma/\gamma$ equal to 100 and 1 is shown in Fig. 8.

Fig.8 (deleted)

Fig. 9 shows the time history of control force in actuator 1 with the upper bound equal to 0.15 for $\sigma/\gamma = 100$. The upper bound is enforced on all the actuators.

Fig.9 8 (deleted)

The changes in the displacement norm with time for 8/y equal to 100 are shown in Fig. 10 for the case of control bound equal to 0.15 and without bound.

Fig.10 (deleted)

Fig. 11 shows the total time required to reduce the displacement norm to 0.02 for three values of $\sigma/\gamma$ and four values of control bound. As the control bound decreases more time is needed to reduce the displacement norm to 0.02 for a given value of $\sigma/\gamma$. The maximum root mean square response for different cases is shown in Fig. 12.

Fig.11 (deleted) Fig.12 (deleted)

## 2. Solution of general linear optimal problem for one control

Now consider the general optimal linear regulator problem with an objective function of minimum time and one control parameter.

Problem Statement. The system is described by a linear differential equation in vector form as,

$$\dot{x} = Ax + Lu \qquad (2\text{-}1)$$

where $x = (x_1, x_2, ..., x_n)$ is the $n$-dimensional state vector , $A = \|a_{ij}\|$ a $n$-dimensional square matrix of constant coefficients, $L$ a column vector which contains $l_1, l_2, ..., l_n$; $u$ a limited control, $|u| \le \zeta$, $\zeta>0$; $x(0) = x_o$, $x(t_k) = x_k$ the initial and final condition, $T = t_k$ represent the end time of process, $t_o = 0$.

It is known the control can have only boundary value in linear system and, if eigenvalues of matrix $A$ is real numbers, the system has only maximum $n-1$ switches [7].

**Problem solution**. The characteristic equation is $|A - \lambda E| = 0$, where $E$ is an unit $n$-dimensional matrix, $\lambda$ is eigenvalues of matrix $A$.

Case A. All eigenvalues λ are real, different, and not equal zero. Using

$$y_i = \sum_{j=1}^{n} e_{ij} x_j \qquad (i=1,2,\ldots,n), \qquad e_{ij} = const_{ij}$$

can convert the equations (2-1) to canonical form

$$\dot{y}_1 = \lambda_1 y_1 + b_1 u; \quad \dot{y}_2 = \lambda_2 y_2 + b_2 u; \quad \ldots \quad \dot{y}_n = \lambda_n y_n + b_n u; \qquad (2\text{-}2)$$

with boundary conditions $y_i(0) = y_{i0}$; $y_i(t_k) = y_{ik}$ .

The optimal control $u = \pm\zeta$ is constant everywhere. If a new variable $z_i = \lambda_i y_i + b_i u$ is introduced, it is possible to write equation (2-2) in form

$$\dot{z}_1 = \lambda_1 z_1; \quad \dot{z}_2 = \lambda_2 z_2; \quad \ldots; \quad \dot{z}_n = \lambda_n z_n; \qquad (2\text{-}3)$$

A solution of equation (2-3) is

$$z_i = \bar{c}_i e^{\lambda_i t} \quad (i = 1, \ldots ,n).$$

Returning to the variable $y$ we can write

$$y_i = c_i e^{\lambda_i t} - b_i u / \lambda_i \quad (i = 1,\ldots,n); \qquad c_i = \bar{c}_i / \lambda_i; \quad \lambda_i \ne 0. \qquad (2\text{-}4)$$

Consider the value $y_1$. The moment when a control parameter is changed it is marked an index "$i$" below and right and left from point $t_i$ by plus and minus sign on top of magnitudes.

Let us suggest, that the control has $k-1$ switches. From continuous condition we have $y_i^- = y_i^+$. Therefore we have

$$c_{ji}^- e^{\lambda_{ji} t_i} - b_i u_i^- / \lambda_i = c_{ji}^+ e^{\lambda_{ji} t_i} - b_i u_i^+ / \lambda_i \qquad (i = 1,\ldots,n) . \qquad (2\text{-}5)$$

From (2-5)

$$c_{ji}^+ = c_{ji}^- + e^{-\lambda_j t_i}(u_i^+ - u_i^-) b_j / \lambda_i \qquad (i = 1,...,k\text{-}1). \qquad (2\text{-}6)$$

The value $c_{j,i+1}^- = c_{j,i}^+$. From (2-6), we get

$$c_{jk}^- = c_{j0}^+ + \sum_{i=1}^{k-1} \frac{b_j}{\lambda_i} e^{-\lambda_i t_i}(u_i^+ - u_i^-) \; . \qquad (2\text{-}7)$$

From the first equation (2-4) and boundary conditions for $y_i$, we find

$$c_{jk}^- = e^{-\lambda_j t_k t_k}(y_{jk} + \frac{b_j}{\lambda_j} u_k^-); \qquad c_{j0}^+ = y_{j0} + \frac{b_j}{\lambda_i} u_0 , \qquad (2\text{-}8)$$

where $u_0 = u_0^+$.

Substituting (2-8) to (2-7) we obtain,

$$\sum_{i=1}^{k-1}\frac{b_j}{\lambda_j}e^{-\lambda_i t_i}(u_i^+ - u_i^-) - e^{-\lambda_j t_k t_k}(y_{jk} + \frac{b_j}{\lambda_j}u_k^-) = -y_{j0} + \frac{b_j}{\lambda_i}u_0 \qquad (j=1,\ldots,n) \qquad (2\text{-}9)$$

These equation (2-9) satisfies for all $y_i$ $i$=1,…,$n$.

If to divide the right and left parts of equation (2-9) by $(-2b_j u_o/\lambda_I)$, we find,

$$e^{-\lambda_1 t_1} - e^{-\lambda_1 t_2} + \ldots - (-1)^k e^{-\lambda_1 t_{k-1}} + \frac{1}{2}\left[(-1)^{k-1} + y_{1k}\lambda_1 / b_1 u_0\right]e^{-\lambda_1 t_k} = \frac{1}{2}\left[1 + y_{10}\lambda_i / b_1 u_0\right];$$

$$e^{-\lambda_2 t_1} - e^{-\lambda_2 t_2} + \ldots - (-1)^k e^{-\lambda_2 t_{k-1}} + \frac{1}{2}\left[(-1)^{k-1} + y_{2k}\lambda_2 / b_2 u_0\right]e^{-\lambda_2 t_k} = \frac{1}{2}\left[1 + y_{20}\lambda_2 / b_2 u_0\right];$$

. . . . . . . . . . . . . . . . . . . . . . . . . . . . . (2-10)

$$e^{-\lambda_n t_1} - e^{-\lambda_n t_1} + \ldots - (-1)^k e^{-\lambda_n t_{k-1}} + \frac{1}{2}\left[(-1)^{k-1} + y_{nk}\lambda_n / b_n u_0\right]e^{-\lambda_n t_k} = \frac{1}{2}\left[1 + y_{n0}\lambda_n / b_n u_0\right],$$

where $k=n$.

Noting that $e^{-t_i} = w_i$ , equations (2-10) can be written as

$$w_1^{\lambda_1} - w_2^{\lambda_1} + \ldots - (-1)^k w_{k-1}^{\lambda_1} + \frac{1}{2}\left[(-1)^{k-1} + y_{1k}\lambda_1 / b_1 u_0\right]w_k^{\lambda_1} = \frac{1}{2}\left[1 + y_{10}\lambda_1 / b_1 u_0\right],$$

$$w_1^{\lambda_2} - w_2^{\lambda_2} + \ldots - (-1)^k w_{k-1}^{\lambda_2} + \frac{1}{2}\left[(-1)^{k-1} + y_{2k}\lambda_2 / b_2 u_0\right]w_k^{\lambda_2} = \frac{1}{2}\left[1 + y_{20}\lambda_2 / b_2 u_0\right],$$

. . . . . . . . . . . . . . . . . . . . . . . . (2-11)

$$w_1^{\lambda_n} - w_2^{\lambda_n} + \ldots - (-1)^k w_{k-1}^{\lambda_n} + \frac{1}{2}\left[(-1)^{k-1} + y_{nk}\lambda_n / b_n u_0\right]w_2^{\lambda_n} = \frac{1}{2}\left[1 + y_{n0}\lambda_n / b_n u_0\right].$$

Equation (2-11) is solved in order to find $w_i = w_i(y_o)$. Returning to the original variable $x$, we can write

$$t_i = -\ln w_i(x) .$$

where $x_o$ represents the initial point $x$.

Equations (2-11) are a set of algebraic equations. From boundary conditions we know that

$$t_k \ge t_{k-1} \ge \ldots \ge t_1 > 0 .$$

That implies that

$$0 < w_k \le w_{k-1} \le \ldots \le w_1 \le 1 .$$

For control $u = \pm\xi$ . This implies that equations (2-11) must be solved twice. If $x(t_k) = 0$ (this means $y(t_k) = 0$), the second solution is symmetric about the origin.

The solution of equation (2-11) is easier to evaluate then the classical optimal control solution. In classical theory a researcher must solve a boundary problem for a set of given differential equations and also find a set of unknown Lagrange multipliers. In using equation (2-10) the researcher first establish the required time increments based upon knowledge of the physical situation.

To find the switch surfaces, for $t_1 = 0$, implies thus the trajectory is located on the first switch surface. In this case in equation (2-11) $e^{-t_1} = w_1 = 1$. We then set about solving the first $n$-1 equations (2-11) for $w_2, w_3, \ldots, w_k$ and substitute these solutions into the last equation. This leads to an equation $\Phi_1(y) = 0$. By substituting for $y$ we can find $N_1(x)$. This is the first (n-1)-dimensional switch surface.

Next by substituting $w_1 = 1$, $w_2 = 1$ in the first $n$-1 equations, and solving the first $n$-2 equations for $w_3, w_4, \ldots, w_n$ one can obtain solutions and substitute them into last equation. We thus can find a hyper surface $N_2(x) = 0$. The intersection of this hyper surface with $N_1(x)$ creates the second ($n$-2)-dimensional switch surface. Other switch surfaces can be found in as similar way.

Such an optimal control result can be easily found. In selecting $u_o$, when the state point reaches the switch surface $N_1(x) = 0$, $u_1 = -u_o$. When the state point reaches the switch surface $N_2(x) = 0$, $u_2 = -u_1$ and so on.

If time is deleted from any two of equation (2-4), we obtain a projection of the trajectory on the surface $y_i y_j$

$$y_i = \bar{c}_i (y_i + b_i u/\lambda_i)^{\lambda_i/\lambda_j} - b_i u/\lambda_i, \qquad \bar{c}_i = c_i / c_j^{\lambda_i/\lambda_j} .$$

From (2-2) we can find the boundaries of instability, for positive eigenvalues. For example, if $\lambda_I > 0$, $b_i > 0$, then $y_i(t_k) = 0$. The necessary and sufficient condition unstable solution is given by

$$y_i > \frac{b_i}{\lambda_i}\xi \ ; \qquad y_i < -\frac{b_i}{\lambda_i}\xi .$$

We have only considered cases when the eigenvalues are real, different, and non-equal to zero. Additional cases have been considered in reference [6] .

**Example**. Taking any two of equations (2-1) with eigenvalues $\lambda_1, \lambda_2$ ($\lambda_1 \neq \lambda_2, \lambda_1 \neq 0, \lambda_2 \neq 0, \lambda_1 < \lambda_2$ ), $x(t_k) = 0$, a canonical form of the equations can be expressed as,

$$\dot{y}_1 = \lambda_1 y_1 + u; \qquad \dot{y}_2 = \lambda_2 y_2 + u; \qquad y(0) = y; \quad y(t_k) = 0; \quad |u| \le 1 . \qquad \text{(E1)}$$

Equation (2-11) for $u_o = 1$ can be written as

$$w_1^{\lambda_1} - \frac{1}{2} w_2^{\lambda_1} = \frac{1}{2}\lambda_1 y_1 + \frac{1}{2} \ ; \qquad w_1^{\lambda_2} - \frac{1}{2} w_2^{\lambda_2} = \frac{1}{2}\lambda_2 y_2 + \frac{1}{2} . \qquad \text{(E2)}$$

For $w_1 = 1$ (simplifying for the case $t_1 = 0$), $w_2$ can be obtained from (E2) and thus

$$y_2 - \frac{1}{\lambda_2}\left[1 - (1 - \lambda_1 y_1)^{\lambda_2/\lambda_1}\right] = 0 . \qquad \text{(E3)}$$

For $u_o = -1$ equations (2-11) are

$$w_1^{\lambda_1} - \frac{1}{2} w_2^{\lambda_1} = -\frac{1}{2}\lambda_1 y_1 + \frac{1}{2} \ ; \qquad w_1^{\lambda_2} - \frac{1}{2} w_2^{\lambda_2} = -\frac{1}{2}\lambda_2 y_2 + \frac{1}{2} . \qquad \text{(E4)}$$

Taking $w_1 = 1$, and using $w_2$ from (E2), we find,

$$y_2 + \frac{1}{\lambda_2}\left[1 - (1 + \lambda_1 y_1)^{\lambda_2/\lambda_1}\right] = 0 . \qquad \text{(E5)}$$

Using a continuity condition $y_1(t_k) = y_2(t_k)$, the relations (E3), (E5) can be written as one relation

$$y_2 - \frac{sign\, y_1}{\lambda_2}\left[1 - (1 - \lambda_1 y_1 sign\ y_1)^{\lambda_1/\lambda_2}\right] = 0 . \qquad \text{(E6)}$$

If $y_1 = 0$, $y_2 > 0$, then the relation (E6) is greater then zero. From (E1) we see: $y_2$ will be decrease faster if $u = -1$ for $y_2 > 0$ and $u = +1$ for $y_2 < 0$. This implied that

$$u = -sign\, \lambda_2\ sign\left\{ y_2 - \frac{sign\ y_1}{\lambda_2}\left[1 - (1 - \lambda_1 y_1 sign\, y_1)^{\lambda_1/\lambda_2}\right]\right\} . \qquad \text{(E7)}$$

To find the equations for optimal trajectories. Referring equations (2-4),(2-11) we find

$$y_1 = c_1 e^{\lambda_1 t} - u/\lambda_1 \ ; \qquad y_2 = c_2 e^{\lambda_2 t} - u/\lambda_2 \ ; \qquad y_2 = c\ (y_1 + u/\lambda_1)^{\lambda_2/\lambda_1} - u/\lambda_2 . \qquad \text{(E8)}$$

The last equation in (E8) gives information in the trajectories as shown in figure 13.

Fig13 (deleted)

These trajectories depend upon the signs of $\lambda_1, \lambda_2$. For then $\lambda_1 > 0$, $\lambda_2 < 0$ the non-stability region is $|y_1| > \zeta/\lambda_1$ . For $\lambda_1 < 0$, $\lambda_2 > 0$ the non-stability region is $|y_2| > \zeta/\lambda_2$ .

In fig.14 also shown optimal trajectories. Once again they depend up on the signs of $\lambda_1, \lambda_2$. Returning to the variables $x$, the picture 14 is affined deformity.

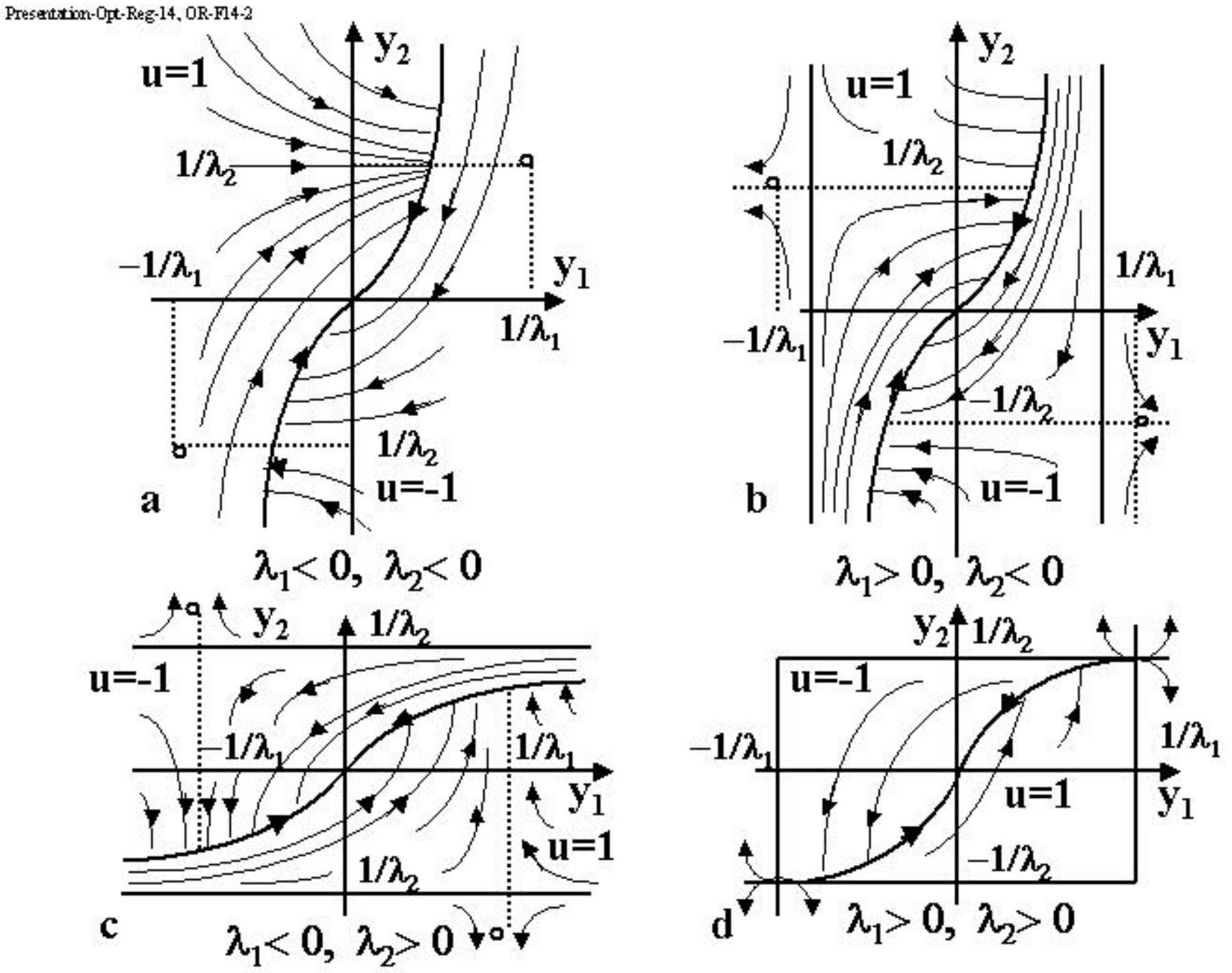

Fig.14.

The offered method allows capture of optimal control.
The other works related to this topics are in [11]-[82].

**Summary**

Two optimum control design methods for suppression of structural vibration having bounded constrains have been compared. The minimum time and quadratic performance index have been used as objective functions. The second approach leads to use of the LQR methodology with bounded control. The introduction of a minimum time controller can be used when the number of actuators equals the number of structural degrees of freedom used in the design of the control system. When the number of actuators is less than the number of degrees of freedom, the minimum time controller becomes mathematically complicated and has been found to be difficult to solve due to the presence of local minimums. The minimum quadratic function controller, with bounded control, can be designed with a fewer number of actuators. A SISO structural control design problem has been solved using both approaches for comparison of trajectories and the time needed to suppress vibrations. The influence of control limitations and the weight coefficient $\sigma/\gamma$ of the structure have been studied. Results indicate that an optimal selection of the weight coefficient $\sigma/\gamma$ can decrease the suppression time up to 2-4 times.

**Recommendations**

If possible, the researcher should try to design the controller for minimum time. If it is very difficult, he can design LQR controller. However in this case the researcher must:

1. Consider limits on the maximum value of the control force.
2. Find the optimal ratio $\sigma/\gamma$ of the weight coefficients.
3. Solve (numerically) at least one time the real (minimum time) problem and compare what may be luck is loss from changing the $T_{min}$ problem to the LQR problem.

# References

(The reader finds some of these works in http://Bolonkin.narod.ru/p65.htm and http://arxiv.org)

1. A.A. Bolonkin, Solution general linear optimal problem with one control. Journal "Pricladnaya Mechanica", v.4, #4, 1968, pgs. 111-122. Moscow (in Russian).
2. **Book**: A.A. Bolonkin, **New methods of optimization and their application**, Moscow, Technical University named Bauman, 1972, pgs.220 (in Russian).
3. A.A. Bolonkin, "Special Extrema in Optimal Control Problems", Akademiya Nauk, Izvestia, Theknicheskaya Kibernetika, No.2, March-April, 1969, pp.187-198. See also English translation in Eng.Cybernetics, n.2, March-April,1969, pp.170-183.
4. A.A. Bolonkin, "A New Approach to Finding a Global Optimum", "New Americans Collected Scientific Reports",Vol.1, 1991. the Bnai Zion Sa Scientists Division, New York.
5. A.A. Bolonkin, N.S. Khot, Optimal Structural Control Design, IAF-94-1.4.206, 45th Congress of the International Astronautical Federation, World Space Congress-1994. October 9-14,1994/Jerusalem, Israel.
6. V. Boltyanski, A. Poznyak, Linear multi-model time-optimization, Journal "Optimal Control Applications and Methods",Vol.23, Issue 3, 2002, pp.141-161.
7. Yunying Mao, Zeyi Liu, The optimal feedback control of the linear-quadratic control problem with a control inequality constraint, Journal "Optimal Control Applications and Methods",Vol.22, Issue 2, 2001, pp.95-109.
8. Heping Hua, Numerical solution of optimal control problems, Journal "Optimal Control Applications and Methods",Vol.21, Issue 5, 2000, pp.233-241.
9. H.Singh, R.H. Broun, D.S. Naidu, Unified approach to linear quadratic regulator with time-scale property, Journal "Optimal Control Applications and Methods",Vol.22, Issue 1, 2001, pp 1-16..
10. B.J. Driessen, N. Sadegh, Minimum-time control of systems with Coulomb friction: near global via mixed integer linear programming, Journal "Optimal Control Applications and Methods",Vol.22, Issue 2, 2001, pp.51-62.
11. **Book:** A.A. Bolonkin, "**Non-Rocket Space Launch and Flight**", Elsevier, London, 2006, 488 ps.
12. A.A. Bolonkin, Optimal trajectory of air and space vehicles, AEAT, #2, pp.193-214, 2004.
13. A.A. Bolonkin, R. Sierakowski, Design of Optimal Regulators, 2nd AIAA "Unmanned Unlimited" Systems, Technologies, and Operations – Aerospace, Land, and sea Conference and Workshop & Exhibit, San Diego, California, 15-18 Sep 2003, AIAA-2003-6638.
14. A.A. Bolonkin, N.S. Khot, Design of Structure control System using Bounded LQG,**Eng.Opt**.,1997. Vol.29.pp, 347-358.
15. A.A. Bolonkin, N.S. Khot, Optimal Bounded Control Design for Vibration Suppression. **Acta Astronautics**, Vol.38, No. 10, pp803-813, 1996.
16. A.A., Bolonkin, Minimum Weight of Control Devices with Bounded LQG Control. **The World Space Congress -96**, June 1-6, 1996, Albuquerque, MN, USA.
17. A.A. Bolonkin, N.S. Khot, Design of Smart Structures with Bounded Controls**, Smart Structures and Materials**, Feb. 25-29,1996,San-Diego, CA.
18. A.A. Bolonkin, Optimum Structural Vibration Control with Bounds on Control Forces, 1995 ASME Design Technical Conference, **15th Biennial Conference on Vibration and Noise,** September 17-21,1995, Boston, MA, USA.
19. A.A. Bolonkin, Design and Optimal Control in Smart Structures. **Conference "Mathematics and Control in Smart Structures**", 26 Feb.-3 March 1995, San Diego, CA, USA.
20. A.A. Bolonkin, Optimal Structural and Control Design. **45th International Actronautical Congress**. Jerusalem, Israel. October 9-14, 1994, IAF-94-I.4.206.

21. A.A. Bolonkin, Optimization of Trajectories of Multistage Rockets. **Investigations of Flight Dynamics**. Moscow, 1965, p. 20-78 (Russian). International Aerospace Abstract A66-23337# (English).
22. A.A. Bolonkin, Special extreme in optimal control. **Akademia Nauk USSR, Izvestiya.** Tekhnicheskaya Kibernetika, No 2, Mar-Apr.,1969, p.187-198. See also *English translation in* **Engineering Cybernetics**, # 2, Mar- Apr.1969, p.170-183, (English).
23. A.A. Bolonkin, Methods of solution for boundary-value problems of Optimal Control Theory. *Translated from* **Prikladnaya Mekhanika**, Vol. 7, No 6,1971, p.639-650, (English).
24. A.A. Bolonkin, Solution of discrete problems of optimal control on the basis of a general minimum principle (Russian. English summary). *Vycisl. Prikl. Mat. (Kiev) Vyp. 7 (1969)*, 121-132. Mathematical Review 771.
25. A.A. Bolonkin, Impulse solution in control problems. (Russian). *Izv. Sibirsk. Otdel. Akad. Nayk USSR*, 1968, No. 13, 63-68. M.R. 7568.
26. A.A. Bolonkin, A certain method of solving optimal problems. *Izv. Sibirsk. Otdel. Akad. Nauk SSSR*. 1970, no.8, p.86-92. M.R. #6163.
27. A.A. Bolonkin, A method for the solution of optimal problems (Russian). *Complex Systems Control*, pp.34-67. Naukova Dumka, Kiev, 1965. M.R. #5535.
28. A.A. Bolonkin, A certain approach to the solution of optimal problems. (Russian. English summary). *Vycisl. Prikl. Mat*. (Kiev). Vyp. 12 (1970), 123-133. M.R. #7940.
29. A.A. Bolonkin, Optimization of parameters of variation problems (Ukrainian. Russian and English summaries). *Dopoviti Akad. Nauk Ukrain. RSR,* 1964, 580-582. M.R. #6352.
30. A.A. Bolonkin, The calculus of variations and a functional equation of Bellman, and an interpretation of Langrange's undetermined multipliers. (Ukrainian. Russian and English summaries). *Dopovidi Akad. Nauk Ukrain., RSR* 1964, 1290-1293. M.R. #5136.
31. A.A. Bolonkin, The extension principle and the Jacobi condition of the variation calculus. (Ukrainian. Russian and English summaries). *Dopovidi Akad. Nauk Ukrain. RSR* 1964. 849-853. M.R. #5117.
32. A.A. Bolonkin, Electrostatic Solar Wind Propulsion System, AIAA-2005-3653. 41 Propulsion Conference, 10-12 July, 2005, Tucson, Arizona, USA.
33. A.A. Bolonkin, Electrostatic Utilization of Asteroids for Space Flight, AIAA-2005-4032. 41 Propulsion Conference, 10-12 July, 2005, Tucson, Arizona, USA.
34. A.A. Bolonkin, Kinetic Anti-Gravitator, AIAA-2005-4504. 41 Propulsion Conference, 10-12 July, 2005, Tucson, Arizona, USA.
35. A.A. Bolonkin, Sling Rotary Space Launcher, AIAA-2005-4035. 41 Propulsion Conference, 10-12 July, 2005, Tucson, Arizona, USA.
36. A.A. Bolonkin, Radioisotope Space Sail and Electric Generator, AIAA-2005-4225. 41 Propulsion Conference, 10-12 July, 2005, Tucson, Arizona, USA.
37. A.A. Bolonkin, Guided Solar Sail and Electric Generator, AIAA-2005-3857. 41 Propulsion Conference, 10-12 July, 2005, Tucson, Arizona, USA.
38. A.A. Bolonkin, Electrostatic Levitation and Artificial Gravity, AIAA-2005-4465. 41 Propulsion Conference, 10-12 July, 2005, Tucson, Arizona, USA.
39. A.A. Bolonkin, Light Multi-reflex Engine, *JBIS*, Vol.57, #9-10, 2004.
40. A.A. Bolonkin, Multi-reflex Propulsion Systems for Space and Air Vehicles and Energy Transfer for Long Distance, JBIS Vol. 57, pp.379-390, 2004.
41. A.A. Bolonkin, High efficiency transfer of mechanical energy. International Energy Conversion Engineering Conference at Providence, RI, Aug.16-19, 2004. AIAA-2004-5660.
42. A.A. Bolonkin, Utilization of Wind Energy at High Altitude. International Energy Conversion Engineering Conference at Providence, RI, Aug.16-19, 2004. AIAA-2004-5756, AIAA-2004-5705.

43. A.A. Bolonkin, Air Cable Transport System, *Journal of Aircraft*, Vol.40, No.2, March-April 2003.
44. A.A. Bolonkin, Earth Accelerators for Space Ships and Missiles, *Journal JBIS*, Vol.56, pp.394-404, 2003.
45. A.A. Bolonkin, Non-Rocket Transportation System for Space Travel, *Journal JBIS*, Vol.56, pp.231-249, 2003.
46. A.A. Bolonkin, Asteroids as propulsion Systems of Space Ships, *JBIS*, Vol.56, pp.98-107, 2003.
47. A.A. Bolonkin, Optimal Inflatable Space Towers with 3-100 km Height, *JBIS*, Vol.56, pp.87-97, 2003.
48. A.A. Bolonkin, Centrifugal Keeper for Space Stations and Satellites. *JBIS*, Vol.56,2003.
49. A.A. Bolonkin, Hypersonic Space Launcher of High Capability, *Actual problems of aviation and space system*. No.1(15), vol.8, pp.45-58, 2003.
50. A.A. Bolonkin, J. Cloutier, Search for Enemy Targets, Technical. Report AFRL-MN-EG-TR-2003-1716, June 2002. 49 p.
51. A.A. Bolonkin, J.Cloutier, Search, Observation, and Attack Problems, *Technical Report* AFRL-MN-EG-TR-2003-1717, 21p. 2002.
52. A.A. Bolonkin, G.Gilyard, Optimal Pitch Thrust Vector Angle and Benefit for all Flight Regimes, NASA/TM- 2000-209021.
53. A.A. Bolonkin, G.Gilyard, Estimated Benefits of Variable-Geometry Wing for Transport Aircraft,NASA/TM-1999-206586.
54. A.A. Bolonkin, Non-Rocket Earth-Moon Transport System, Journal "*Advance in Space Research*", Vol.31/11, pp.2485-2490.
55. A.A. Bolonkin, A High Efficiency Fuselage propeller ("Fusefan") for Subsonic Aircraft, !999 World Aviation Congress, AIAA, #1999-01-5569.
56. A.A. Bolonkin, Air Cable Transport and Bridges, TN 7567, International Air & Space Symposium – The Next 100 Years, 14-17 July 2003, Dayton, Ohio, USA
57. A.A. Bolonkin, "Non-Rocket Space Rope Launcher for People", IAC-02-V.P.06, 53rd International Astronautical Congress. The World Space Congress – 2002, 10-19 Oct 2002/Houston, Texas, USA.
58. A.A. Bolonkin, "Non-Rocket Missile Rope Launcher", IAC-02-IAA.S.P.14, 53rd International Astronautical Congress. The World Space Congress – 2002, 10-19 Oct 2002/Houston, Texas, USA.
59. A.A. Bolonkin, "Inexpensive Cable Space Launcher of High Capability", IAC-02-V.P.07, 53rd International Astronautical Congress. The World Space Congress – 2002, 10-19 Oct. 2002/Houston, Texas, USA.
60. A.A. Bolonkin, "Hypersonic Launch System of Capability up 500 tons per day and Delivery Cost $1 per Lb". IAC-02-S.P.15, 53rd International Astronautical Congress. The World Space Congress – 2002, 10-19 Oct 2002/Houston, Texas, USA.
61. A.A. Bolonkin, "Employment Asteroids for Movement of Space Ship and Probes". IAC-02-S.6.04, 53rd International Astronautical Congress. The World Space Congress – 2002, 10-19 Oct. 2002/Houston, Texas, USA.
62. A.A. Bolonkin, "Optimal Inflatable Space Towers of High Height". COSPAR-02 C1.1-0035-02, 34th Scientific Assembly of the Committee on Space Research (COSPAR). The World Space Congress – 2002, 10-19 Oct 2002/Houston, Texas, USA.
63. A.A. Bolonkin, "Non-Rocket Earth-Moon Transport System", COSPAR-02 B0.3-F3.3-0032-02, 02-A-02226, 34th Scientific Assembly of the Committee on Space Research (COSPAR). The World Space Congress – 2002, 10-19 Oct 2002/Houston, Texas, USA.

64. A.A. Bolonkin, "Non-Rocket Earth-Mars Transport System", COSPAR-02 B0.4-C3.4-0036-02, 34th Scientific Assembly of the Committee on Space Research (COSPAR). The World Space Congress – 2002, 10-19 Oct 2002/Houston, Texas, USA.
65. A.A. Bolonkin, "Transport System for delivery Tourists at Altitude 140 km". IAC-02-IAA.1.3.03, 53rd International Astronautical Congress. The World Space Congress – 2002, 10-19 Oct. 2002/Houston, Texas, USA.
66. A.A. Bolonkin, Space Cable Launchers, TN 8057, International Air & Space Symposium – The Next 100 Years, 14-17 July 2003, Dayton, Ohio, USA.
67. A.A. Bolonkin, Geometry-Based Feasibility Constraints for Single Pursuer Multiple Evader Problems, 2nd AIAA "Unmanned Unlimited" Systems, Technologies, and Operations – Aerospace, Land, and sea Conference and Workshop & Exhibit, San Diego, California, 15-18 Sep 2003, AIAA-2003-6638.
68. A.A. Bolonkin, The Simplest Space Electric Generator and Motor with Control Energy and Thrust, **45th International Astronautical Congress,** Jerusalem. Israel. Oct. 9-14, 1994, IAF-94-R.1.368
69. A.A. Bolonkin, Space Electric Generator, run by Solar Wing. **The World Space Congress**, Washington, DC, USA, 28 Aug.-5 Sept., 1992, IAF-92-0604.
70. A.A. Bolonkin, Simple Space Nuclear Reactor Motors and Electric Generators Running on Radioactive Substances, **The World Space Congress,** Washington, DC, USA, 28 Aug.- 5 Sept., 1992, IAF-92-0573.
71. A.A. Bolonkin, A Space Motor Using Solar Wind Energy (Magnetic Particle Sail). **The World Space Congress**, Washington, DC, USA, 28 Aug.- 5 Sept., 1992, IAF-0615.
72. A.A. Bolonkin, Aviation, Motor and Space Designs. **Emerging Technology in the Soviet Union,** 1990, Delphic Ass., Inc., pp.32-80.
73. **Book:** A.A. Bolonkin, **The Development of Soviet Rocket Engines,** 1991,Delphic Ass. Inc.,122 p., Washington,(Eng).
74. A.A. Bolonkin, **New Way of Generation of Electrical Energy in Space**. Report ESTI, 1988, 109p. (Soviet Classified Project).
75. A.A. Bolonkin, **New Rotor Internal Combustion Engine**. Report ESTI,1988,75p. (Soviet Classif. project),(Rus.)
76. A.A. Bolonkin, **Supersonic VTOL fighter-helicopter**. Report ESTI, 1988, 120p. (Soviet Clas. Project),
77. A.A. Bolonkin, **New ground-effect vehicle**. Report ESTI, 1987, 85p. (Soviet Classified Project),
78. A.A. Bolonkin, **New rotor hydraulic transmission for cars and other machines**. Report ESTI, 1986. (Soviet Classified Project).
79. A.A. Bolonkin, Investigation of the take off dynamics of a VTOL aircraft. **Investigations of Flight Dynamics**. Moscow, 1965, p. 119-147 ( Russian). International Aerospace Abstract A66-23339# (English).
80. A.A. Bolonkin, Theory of lifting body with controllable radial force. **Investigations of Flight Dynamics.** Moscow**, 1965, p.79-118,** (Russian).International Aerospace Abstract A66-23338#(Eng).
81. **Book:** A.A. Bolonkin, **Theory of Flight Models**, Moscow, Association of Army, Air Force, and NAVY, 328p. 1962 (in Russian)

## Nomenclature

$A$ is the $2n\times 2n$ plant matrix in liner problem
$a_{ij}$ is members of matrix A

$B$ is $2n \times p$ control matrix in liner problem
$b_{jk}$ is members of matrix B
$C$ is $q \times 2n$ out matrix
$C_i$ are weight coefficients
$c$ is constant
$F_i$ is the magnitude of the bounds for each controller
$F_0$ is function of initial conditions
$f$ is the control force vector of dimension $p$
$|f| \leq 1$ is the control in linear problem
$H$ is Hamiltonian
$I$ is the functional (objective function),
$\overline{M}, \overline{E}, \overline{K}$ are diagonal square matrices
$P$ is a $2n$-dimensional unknown matrix
$\underline{Q}$ is state weighting matrices
$\underline{R}$ is control weighting matrices
$R_x(t)$ is norm of displacement
$T$ is final time
$t$ is time (variable)
$t_1$, $t_2$, - boundary condition
$u$ is the vector defines the structural response.
$v$ is a $p$-dimensional vector of control forces
$x$ is a $n$-dimensional vector of state in general problem,
$x$ is the state vector of dimension $2n$ in linear problem.
$x(0)$ is the initial state vector
$x(T)$ is the final state of the system.
$x_1$, $x_2$ - boundary condition
$\zeta$ is the vector of modal damping factors
$\lambda(t)$ is a $n$-dimensional vector unknown coefficient
$\lambda$ is eigenvalues of matrix $A$
$\psi = \psi(t,x)$ is special function
$\omega$ is the vector of structural frequencies